\documentclass[a4paper,11pt]{article}
\usepackage{latexsym, amssymb, amsfonts,amsmath}
\usepackage{amssymb,latexsym}
\usepackage{amsfonts}
\usepackage{dsfont}
\usepackage{color}

\usepackage{amsmath,amsthm,graphicx}
\usepackage[all]{xy}
\newcommand{\noise}{\mathsf{k}}

\newcommand{\bea}{\begin{eqnarray}}
\newcommand{\eea}{\end{eqnarray}}
\newcommand{\brray}{\begin{array}}
\newcommand{\erray}{\end{array}}
\newcommand{\benu}{\begin{enumerate}}
\newcommand{\eenu}{\end{enumerate}}

\newcommand{\eprf}{\end{proof}}
\newcommand{\newsection}[1]{\setcounter{equation}{0} \setcounter{dfn}{0}

\begin{document}
asdasd 
\end{document}

\section{#1}}
\newcommand{\be}{\begin{equation}}
\newcommand{\ee}{\end{equation}}
\newcommand{\CB}{CB}
\newcommand{\cnoise}{\widehat{\noise}}
\newcommand{\id}{\textup{id}}
\newcommand{\I}{\textup{I}}
\newcommand{\tud}{\textup{-}}
\newcommand{\msfhh}{\mathsf{h}}
\newcommand{\cmsfhh}{\widehat{\mathsf{h}}}
\newcommand{\mfs}{\rho}
\newcommand{\cinit}{\widehat{\init}}
\newcommand{\msfh}{\mathsf{H}}
\newcommand{\otul}{\underline{\ot}}
\newcommand{\qiso}{\mathbb{I}_{\mathbb{M}}}
\newcommand{\oqiso}{\mathbb{I}_{\mathbb{M},0}}
\newcommand{\qisored}{\mathbb{I}_{\mathbb{M,red}}}
\newcommand{\manifold}{\mathbb{M}}
\newcommand{\catq}{\textup{{\bf Q}}_\manifold} 
\newcommand{\mor}{\textup{Mor}}
\newcommand{\hopf}{(A_0,\Delta_0)}
\newcommand{\full}{(A,\Delta)}
\newcommand{\otol}{\overline{\ot}}
\newcommand{\qsf}{\mathbb{Q}\mathbb{S}\mathbb{F}}
\newcommand{\qlp}{\mathbb{Q}\mathbb{L}\mathbb{P}}
\newcommand{\ogamma}{\Gamma_0}
\newcommand{\qsc}{\mathbb{Q}\mathbb{S}_u\mathbb{C}}
\newcommand{\ben}{\begin{equation*}\begin{split}}
\newcommand{\een}{\end{equation*}}
\newcommand{\esp}{\end{split}}
\newcommand{\qscc}{\mathbb{Q}\mathbb{S}_c\mathbb{C}}
\newcommand{\msf}{\mathsf{}}
\newcommand{\sro}{\mathsf{F}}
\newcommand{\qg}{\mathbb{G}}
\newcommand{\oqg}{\mathbb{G}_0}
\newcommand{\cyc}{1^{(2)}}
\newcommand{\mult}{\mathsf{W}}
\newcommand{\coproduct}{\Delta}
\DeclareMathOperator{\linear}{Lin}
\newcommand{\ocoproduct}{\Delta_0}
\newcommand{\oinit}{\mathfrak{h}_{00}}
\newcommand{\cbm}{CB^\coproduct_r(A;A\ot_M B)}
\newcommand{\cb}{CB^\coproduct_r(A;A\ot B)}
\newcommand{\lin}{L^\coproduct(A_0;A\ot B)}
\newcommand{\linm}{L^\coproduct(A_0;A\ot_M B)}
\newcommand{\linalg}{L^\coproduct(A_0;A_0\otul B)}
\newcommand{\m}{\mathbb{M}_\mathbb{R}^\Delta(A_r)}
\newcommand{\WW}{{\mathds{V}\!\!\text{\reflectbox{$\mathds{V}$}}}}
\newcommand{\Ww}{\mathds{W}}
\newcommand{\wW}{\text{\reflectbox{$\Ww$}}\:\!}
\newcommand{\rdual}{\widehat{R}_u}
\newcommand{\duh}{\widehat{\IH}}
\newcommand{\braH}{\innerl H|}
\newcommand{\ketH}{|H\innerr}

\newtheorem{dfn}{Definition}[section]
\newtheorem{thm}[dfn]{Theorem}
\newtheorem{nota}[dfn]{Notation}
\newtheorem{lmma}[dfn]{Lemma}
\newtheorem{hypo}[dfn]{Hypothesis}
\newtheorem{ppsn}[dfn]{Proposition}
\newtheorem{crlre}[dfn]{Corollary}
\newtheorem{xmpl}[dfn]{Example}
\newtheorem{rmrk}[dfn]{Remark}
\newcommand{\bdfn}{\begin{dfn}}
\newcommand{\bthm}{\begin{thm}}

\newtheorem{ass}{Assumption}

\newcommand{\qito}{\Delta^{\textup{QS}}}
\newcommand{\init}{\mathfrak{h}}

\newcommand{\blr}{\begin{list}{$($\roman{cnt1}$)$} {\usecounter{cnt1}
        \setlength{\topsep}{0pt} \setlength{\itemsep}{0pt}}}
\newcommand{\bla}{\begin{list}{$($\alph{cnt2}$)$} {\usecounter{cnt2}
       \setlength{\topsep}{0pt} \setlength{\itemsep}{0pt}}}
\newcommand{\bln}{\begin{list}{$($\arabic{cnt3}$)$} {\usecounter{cnt3}
                \setlength{\topsep}{0pt} \setlength{\itemsep}{0pt}}}
\newcommand{\el}{\end{list}}
\newcommand{\no}{\noindent}

\newcommand{\blmma}{\begin{lmma}}
\newcommand{\bppsn}{\begin{ppsn}}
\newcommand{\bcrlre}{\begin{crlre}}
\newcommand{\bxmpl}{\begin{xmpl}}
\newcommand{\brmrk}{\begin{rmrk}}
\newcommand{\edfn}{\end{dfn}}
\newcommand{\ethm}{\end{thm}}
\newcommand{\elmma}{\end{lmma}}
\newcommand{\ptrace}{Tr_{|f><g|}}
\newcommand{\Ptrace}{Tr_{|f^{\ot^m}><g^{\ot^n}|}}
\newcommand{\till}{\widetilde{\cll}}
\newcommand{\eppsn}{\end{ppsn}}
\newcommand{\ecrlre}{\end{crlre}}
\newcommand{\exmpl}{\end{xmpl}}
\newcommand{\ermrk}{\end{rmrk}}

\newcommand{\IA}{{I\! \! A}}
\newcommand{\IB}{{I\! \! B}}
\newcommand{\IC}{\mathbb{C}}

\newcommand{\ID}{{I\! \! D}}
\newcommand{\IE}{\mathbb{E}}
\newcommand{\IF}{{I\! \! F}}
\newcommand{\IG}{{I\! \! G}}
\newcommand{\IH}{\mathbb{H}}
\newcommand{\II}{{I\! \! I}}
\newcommand{\deltaext}{\widetilde{\Delta_u}}
\newcommand{\IK}{{I\! \! K}}
\newcommand{\IL}{\mathbb{L}}
\newcommand{\IM}{{I\! \! M}}
\newcommand{\IN}{{I\! \! N}}
\newcommand{\IO}{{I\! \! O}}
\newcommand{\IP}{{I\! \! P}}
\newcommand{\IQ}{\mathbb{Q}}
\newcommand{\mc}{\mathcal}
\newcommand{\IR}{\mathbb{R}}
\newcommand{\IS}{{I\! \! S}}
\newcommand{\IT}{\mathbb{T}}
\newcommand{\IU}{\mathbb{U}}
\newcommand{\IV}{{I\! \! V}}
\newcommand{\IW}{\mathbb{W}}
\newcommand{\IX}{{I\! \! X}}
\newcommand{\IY}{{I\! \! Y}}
\newcommand{\IZ}{\mathbb{Z}}

\newcommand{\al}{\alpha}
\newcommand{\bta}{\beta}
\newcommand{\gma}{\gamma}
\newcommand{\Dlt}{\Delta}
\newcommand{\dlt}{\delta}
\newcommand{\veps}{\varepsilon}
\newcommand{\eps}{\epsilon}
\newcommand{\kpa}{\kappa}
\newcommand{\lmd}{\lambda}
\newcommand{\Lmd}{\Lambda}
\newcommand{\omg}{\omega}
\newcommand{\innerl}{\langle}
\newcommand{\innerr}{\rangle}
\newcommand{\Omg}{\Omega}
\newcommand{\htpi}{\hat{\pi}}
\newcommand{\sgm}{\sigma}
\newcommand{\Sgm}{\Sigma}
\newcommand{\tta}{\theta}
\newcommand{\Tta}{\Theta}
\newcommand{\vtta}{\vartheta}
\newcommand{\zta}{\zeta}
\newcommand{\del}{\partial}
\newcommand{\Gma}{\Gamma}
\newcommand{\cla}{{\cal A}}
\newcommand{\clb}{{\cal B}}
\newcommand{\clc}{{\cal C}}
\newcommand{\cld}{{\cal D}}
\newcommand{\cle}{{\cal E}}
\newcommand{\clf}{{\cal F}}
\newcommand{\clg}{{\cal G}}
\newcommand{\clh}{{\cal H}}
\newcommand{\cli}{{\cal I}}
\newcommand{\clj}{{\cal J}}
\newcommand{\clk}{{\cal K}}
\newcommand{\cll}{{\cal L}}
\newcommand{\clm}{{\cal M}}
\newcommand{\cln}{{\cal N}}
\newcommand{\clp}{{\cal P}}
\newcommand{\clq}{{\cal Q}}
\newcommand{\clr}{{\cal R}}
\newcommand{\cls}{{\cal S}}
\newcommand{\clt}{{\cal T}}
\newcommand{\clu}{{\cal U}}
\newcommand{\clv}{{\cal V}}
\newcommand{\clw}{{\cal W}}
\newcommand{\clx}{{\cal X}}
\newcommand{\cly}{{\cal Y}}
\newcommand{\clz}{{\cal Z}}
\newcommand{\capA}{\hspace{-.05in}{\bf A}}
\newcommand{\ap}{\mathbb{A}\mathbb{P}}

\def\wA{\widetilde{A}}
\def\wB{\widetilde{B}}
\def\wC{\widetilde{C}}
\def\wD{\widetilde{D}}
\def\wE{\widetilde{E}}
\def\wF{\widetilde{F}}
\def\wG{\widetilde{G}}
\def\wH{\widetilde{H}}
\def\wI{\widetilde{I}}
\def\wJ{\widetilde{J}}
\def\wK{\widetilde{K}}
\def\wL{\widetilde{L}}
\def\wM{\widetilde{M}}
\def\wN{\widetilde{N}}
\def\wO{\widetilde{O}}
\def\wP{\widetilde{P}}
\def\wQ{\widetilde{Q}}
\def\wR{\widetilde{R}}
\def\wS{\widetilde{S}}
\def\wT{\widetilde{T}}
\def\wU{\widetilde{U}}
\def\wV{\widetilde{V}}
\def\wX{\widetilde{X}}
\def\wY{\widetilde{Y}}
\def\wZ{\widetilde{Z}}

\def\what{\widehat}

\def\ah{{\cal A}_h}
\def\a*{{\cal A}_{h,*}}
\def\B{{\cal B}(h)}
\def\B1{{\cal B}_1(h)}
\def\b{{\cal B}^{\rm s.a.}(h)}
\def\b1{{\cal B}^{\rm s.a.}_1(h)}
\def\Ap{{\cal A}^{\perp}}

\newcommand{\itt}{\int \limits}
\newcommand{\qgdual}{\widehat{\qg}}
\newcommand{\codu}{\widehat{\Delta}_u}
\newcommand{\codr}{\widehat{\Delta}_r}
\newcommand{\wht}{\widehat}
\newcommand{\rglr}{\Re}
\newcommand{\ot}{\otimes}
\newcommand{\bgots}{\bigotimes}
\newcommand{\raro}{\rightarrow}
\newcommand{\RARO}{\Rightarrow}
\newcommand{\lraro}{\Longrightarrow}
\newcommand{\sbs}{\subset}
\newcommand{\seq}{\subseteq}
\newcommand{\llra}{\Longleftrightarrow}
\newcommand{\ul}{\underline}
\newcommand{\ol}{\overline}
\newcommand{\lgl}{\langle}
\newcommand{\rgl}{\rangle}
\newcommand{\nn}{\nonumber}
\newcommand{\tnsr}{\mbox{$\bigcirc\hspace{-0.89em}\mbox{\raisebox%
{-.43ex}{$\top$}}\;$}}
\newcommand{\gel}{\gtreqqless}
\newcommand{\leg}{\lesseqqgtr}
\newcommand{\half}{\frac{1}{2}}
\newcommand{\one}{1\!\!1}
\newcommand{\aff}{\mbox{{\boldmath $\eta$}}}
\newcommand{\NI}{\noindent}
\newcommand {\CC}{\centerline}
\def \qed {$\Box$}
\newcommand{\dsp}{\displaystyle}
\newcommand{\vsp}{\vskip 1em}
\newcommand{\cg}{C^u_0(\qg)}
\newcommand{\cgdual}{C^u_0(\qgdual)}
\newcommand{\btil}{\widetilde{B}}
\newcommand{\bg}{B^{(\ast)}(\qg)}
\newcommand{\ww}{\WW^{\prime\ast}}
\newtheorem{theoremA}{Theorem}
\newtheorem{theoremB}{Theorem}
\newtheorem{theoremC}{Theorem}

\begin{document}
\begin{center}
{\large{\bf On a quantum version of Ellis joint continuity theorem}} 
\end{center}

\begin{center}
 Biswarup Das\footnote{Institute of Mathematics of the Polish Academy of Sciences; email: b.das@impan.pl}~and~Colin Mrozinski\footnote{Institute of Mathematics of the Polish Academy of Sciences; email: cmrozinski@impan.pl} 
\end{center}

\begin{abstract}
{\it We give a necessary and sufficient condition on a compact semitopological quantum semigroup which turns it into a compact quantum group. In particular, we obtain a generalisation of Ellis joint continuity theorem. We also investigate the question of the existence of the Haar state on a compact semitopological quantum semigroup and prove a ``noncommutative" version of the converse Haar's theorem.} 
\end{abstract}
\begin{center}
\end{center}
\section{Introduction}\label{InTroduction}
Compact semitopological semigroups \emph{i.e.} compact semigroups with separately continuous product arise naturally in the study of weak almost periodicity in locally compact groups. For example, the weakly almost periodic functions on a locally compact group $G$ form a commutative C*-algebra WAP$(G)$ whose character space $G^{\text{WAP}}$ becomes a compact semitopological semigroup. From an abstract algebraic perspective, one can come up with necessary and sufficient conditions on a semigroup, which make it embeddable (by which we mean an injective group homomorphism) into a group (for example, Ore's Theorem for semigroups in \cite[pp. 35] {ore}). However in general, such abstract conditions do not produce a topological group. The added difficulty in obtaining a topological group from a semitopological semigroup, lies in the fact that not only the semigroup should have a neutral element and existence of inverse of all elements, but one also requires the joint continuity of the product and continuity of the inverse.

In fact, the transition from semitopological semigroups (\emph{i.e.} separate continuity of the product) to topological groups (\emph{i.e.} joint continuity of the product, existence and continuity of the inverse and existence of a neutral element) may be achieved in two different ways: 
\begin{itemize}
 \item[(a)]
 A (locally) compact semitopological semigroup becomes a topological group by requiring that the semigroup is algebraically (\emph{i.e.} as a set) a group. This is known as Ellis joint continuity theorem (see \cite{Ellis}), which plays a fundamental role in the theory of semitopological semigroups.
 \item[(b)]
 A (locally) compact semitopological semigroup with a neutral element and an invariant measure of full support is a (locally) compact group. This is known as converse Haar's theorem (see \cite{mukherjea}).
\end{itemize}

In the recent years, ``noncommutative joint continuity" has been extensively studied under the heading of the ``topological quantum groups", which we shall take to mean C*-bialgebras, probably with additional structures, such as (locally) compact quantum groups. A recent work (see \cite{matt_separate_continuity}) addresses the issue of ``noncommutative separate continuity", through the formulation of weak almost periodicity of Hopf von Neumann algebras, and in particular, it gives a notion of compact semitopological quantum semigroups \cite[Definition 5.3]{matt_separate_continuity} which are quantum analogues of compact semitopological semigroups. 

In this paper, we aim at connecting ``noncommutative separate continuity" with ``noncommmutative joint continuity" through studying noncommutative analogues of (a) and (b).
Our main aim in this paper is to show that a {\it compact semitopological quantum semigroup with ``weak cancellation laws" (see Definition \ref{weak cancellation})} can be looked upon as a quantum analogue of {\it a compact semitopological semigroup which is algebraically a group}. Upon establishing this, we prove a quantum version of Ellis joint continuity theorem:

\renewcommand{\thetheoremA}{\ref{Quantum Ellis theorem}}
\begin{theoremA}[{\it{\bf Quantum Ellis joint continuity theorem}}]
A compact semitopological quantum semigroup is a compact quantum group if and only if it satisfies weak cancellation laws. 
\end{theoremA}

\noindent Specializing this to a compact semitopological semigroup, we get a new proof of the Ellis joint continuity theorem (\cite{Ellis}). It also extends the previously known equivalence of weak cancellation laws and Woronowicz cancellation laws in the context of compact quantum (topological) semigroups \cite[Theorem 3.2]{weak cacellation}. 

We also prove a converse Haar's theorem for compact semitopological quantum semigroup:

\renewcommand{\thetheoremB}{\ref{faithful invariant mean gives CQG}}
\begin{theoremB}[{\it{\bf Quantum converse Haar's theorem}}]
A compact semitopological quantum semigroup with a bounded counit and an invariant state is a compact quantum group.
\end{theoremB}

These two theorems are deduced from the following general result which we prove in Section \ref{main results}:
\renewcommand{\thetheoremC}{\ref{S is a compact quantum group}}
\begin{theoremC}
A compact semitopological quantum semigroup $\mathbb{S}:=(A,\Delta)$ satisfying the assumptions:
\begin{enumerate}
 \item there exists an invariant mean $h\in A^\ast$ on $\mathbb{S}$,
 \item the sets $\linear\{a\star hb:~a,b\in A\}$ and $\linear\{ha\star b:~a,b\in A\}$ are norm dense in $A$
\end{enumerate} 
 is a compact quantum group.   
\end{theoremC}



The paper is organised as follows:  We introduce some terminology and notations in Section 2. Section 3 is devoted to a series of results, leading to the main theorem of this section (Theorem \ref{S is a compact quantum group}). Finally in Sections 4 and 5 we give possible applications of the result obtained in Section 3. In Section 4 we prove a generalization of Ellis joint continuity theorem (Theorem \ref{Quantum Ellis theorem}) and in Section 5 we prove a quantum version of converse Haar's theorem (Theorem \ref{faithful invariant mean gives CQG}).
\section*{Acknowledgement}
Both the authors acknowledge the support of Institute of Mathematics of the Polish Academy of Sciences (IMPAN) and Warsaw Center of Mathematics and Computer Science (WCMS). They would like to thank Dr. Matthew Daws and Dr. Adam Skalski for various mathematical discussions leading to improvements of the results. They would also like to thank Dr. Tomasz Kania for various useful comments at an earlier version of the paper.
\section{Compact semitopological quantum semigroup}
Throughout the paper, the symbols $\otol$ and $\ot$ will denote respectively von Neumann algebraic tensor product (spatial) and minimal C*-algebraic tensor product (spatial). For a Hilbert space $H$, $B_0(H)$ will denote the C*-algebra of compact operators on $H$. By SOT* topology on $B(H)$, we will mean the strong* operator topology. For two Hilbert spaces $H_1$ and $H_2$, $H_1\ot H_2$ will denote the Hilbert space tensor product. For a C*-algebra $A$, we will identify the bidual $A^{\ast\ast}$ with the universal enveloping von Neumann algebra of $A$. For $\omega\in A^\ast$, the symbol $\widetilde{\omega}$ will denote the unique normal extension of $\omega$ to $A^{\ast\ast}$, and for a non-degenerate *-homomorphism $\pi:A\longrightarrow B(H)$, $\widetilde{\pi}$ will denote the normal extension of $\pi$ to $A^{\ast\ast}$. For a von Neumann algebra $N$ with predual $N_\ast$, we will refer to the weak* topology on $N$ as the ultraweak topology on $N$.    

We begin by briefly recalling the definition of a compact semitopological quantum semigroup, building upon the motivations coming from the classical situation. All this has been extensively explained in \cite{matt_separate_continuity}.

\subsection{Classical compact semitopological semigroup}\label{subsection1}
We call a compact semigroup $S$ a semitopological semigroup if the multiplication in $S$ is separately continuous \emph{i.e.} for each $s\in S$, $t\mapsto ts$ and $t\mapsto st$ are continuous. We collect some facts about bounded, separately continuous functions on $S\times S$, which have been discussed in details in Section 3 of \cite{matt_separate_continuity} (also see \cite{matt_separate_continuity1}). 

Let $SC(S\times S)$ denotes the algebra of bounded, separately continuous functions on $S\times S$. Define a map \[\Delta:C(S)\longrightarrow SC(S\times S)~ \text{by}~ \Delta(f)(s,t):=f(st)\quad (\forall~f\in C(S)~\text{and}~ s,t\in S).\]
An argument similar to the proof of Proposition 4.1 in \cite{matt_separate_continuity1} shows that $\Delta$ can be viewed as a unital *-homomorphism from $C(S)$ to $C(S)^{\ast\ast}\otol C(S)^{\ast\ast}$, where $C(S)^{\ast\ast}$ is the bidual of $C(S)$, identified with the universal enveloping von Neumann algebra of $C(S)$. Moreover, it follows from the discussions in Sections 3 and 5 of \cite{matt_separate_continuity} that $(\widetilde{\mu}\ot\id)(\Delta(f))\in C(S)$ and also $(\id\ot\widetilde{\mu})(\Delta(f))\in C(S)$.

Associativity of the product in $S$ implies that $(\widetilde{\Delta}\ot\id)\circ\Delta=(\id\ot\widetilde{\Delta})\circ\Delta$, where $\widetilde{\Delta}$ is the normal extension of $\Delta$ to $C(S)^{\ast\ast}$.

\subsection{Compact semitopological quantum semigroup}
The discussion in Subsection \ref{subsection1} allows us to formulate the following definition of a compact semitopological quantum semigroup (see Definition 5.3 in \cite{matt_separate_continuity}):

\bdfn\label{definition of compact semitopological quantum semigroup}
A compact semitopological quantum semigroup is a pair $\mathbb S := (A,\Delta)$ where 
\begin{itemize}
 \item 
$A$ is a unital C*-algebra, considered as a norm closed C*-subalgebra of $A^{\ast\ast}$.
\item
$\Delta:A\longrightarrow A^{\ast\ast}\otol A^{\ast\ast}$ is a unital $\ast$-homomorphism satisfying 
\[
(\widetilde{\Delta}\ot\id)\circ\Delta=(\id\ot\widetilde{\Delta})\circ\Delta, 
\]
where $\widetilde{\Delta}$ is the normal extension of $\Delta$ to $A^{\ast\ast}$. As usual, we will refer to $\Delta$ as the coproduct of $\mathbb S$.
\item
For $\omega\in A^\ast$, 
\[
(\widetilde{\omega}\ot\id)(\Delta (x))\in A~;~(\id\ot\widetilde{\omega})(\Delta (x))\in A\quad(\forall x\in A), 
\]
where $\widetilde{\omega}$ is the normal extension of $\omega$ to $A^{\ast\ast}$.
\end{itemize}
\edfn
\begin{rmrk}
In the notation of Definition 5.3 in \cite{matt_separate_continuity}, $\Delta:A\longrightarrow A\stackrel{\text{sc}}{\ot}A$. However, in this paper we will not use this notation. 
\end{rmrk}

\begin{xmpl} $\ $
\begin{enumerate}
 \item If $S$ is a compact semitopological semigroup, then $(C(S), \Delta)$ is a compact semitopological quantum semigroup.

 Conversely, if the C*-algebra $A$ in Definition \ref{definition of compact semitopological quantum semigroup} is commutative, then it follows that $A=C(S)$ for some compact semitopological semigroup $S$.
 \item A compact quantum group is a compact semitopological quantum semigroup.
 \item 

The following theorem implies that a unital C*-Eberlein algebra (Definition 3.6 in \cite{das-daws}) is a compact semitopological quantum semigroup (also see the comments after Definition 5.3 in \cite{matt_separate_continuity}). 
\end{enumerate}
\end{xmpl}
\bthm\label{answer to Matt's question}
Suppose $(A,\Delta,V,\clh)$ is a C*-Eberlein algebra (see Definition 3.6 in \cite{das-daws}) with $A$ being unital. Then $(A,\Delta)$ is a compact semitopological quantum semigroup.
\ethm

\begin{proof}
Let $\xi,\eta\in\clh$ and consider the functional $\omega_{\xi,\eta}\in B(\clh)_\ast$. Then we have that 
\[
\Delta((\id\ot\omega_{\xi,\eta})(V))=\sum_{i\in\cli}(\id\ot\omega_{e_i,\eta})(V)\ot(\id\ot\omega_{\xi,e_i})(V), 
\]
for some orthonormal basis $\{e_i\}_{i\in\cli}$ of $\clh$, where the sum obviously converges in the ultraweak topology of $A^{\ast\ast}\otol A^{\ast\ast}$ but also converges in $A^{\ast\ast}\ot_{\text{eh}}A^{\ast\ast}$, the extended Haagerup tensor product (see \cite{bletcher}, where it is called the weak*-Haagerup tensor product). Now by Lemma 2.5 in \cite{bletcher}, it follows that if $f\in A^\ast$, then both the sums
$\sum_i (\id\ot\omega_{e_i,\eta})(V)f((\id\ot\omega_{\xi,e_i})(V))$ and $\sum_i f((\id\ot\omega_{e_i,\eta})(V))(\id\ot\omega_{\xi,e_i})(V)$ converge in the norm topology of $A^{\ast\ast}$. This coupled with the definition of C*-Eberlein algebra (Definition 3.6 in \cite{das-daws}) imply that $(\id\ot f)\big(\Delta((\id\ot\omega_{\xi,\eta})(V))\big)\in A$ and also $(f\ot\id)\big(\Delta((\id\ot\omega_{\xi,\eta})(V))\big)\in A$.

Since the set $\{(\id\ot\omega_{\xi,\eta})(V):~\xi,\eta\in\clh\}$ is norm dense in $A$, it follows that $(\id\ot f)(\Delta(a))\in A$ and $(f\ot\id)(\Delta(a))\in A$ for all $a\in A, f\in A^\ast$. This proves the claim.
\end{proof}

\section{The general framework}\label{main results}
In this section we obtain a general result (Theorem \ref{S is a compact quantum group}) on when a compact semitopological quantum semigroup is a compact quantum group. This will be extensively used in the following sections to obtain generalizations of Ellis joint continuity theorem and converse Haar's theorem. 
\subsection{Some observations}\label{some observations}

Let $\mathbb{S}:=(A,\Delta)$ be a compact semitopological quantum semigroup as in Definition \ref{definition of compact semitopological quantum semigroup}. The coproduct $\Delta$ determines a multiplication in $A^\ast$ given by \[\lambda\star\mu:=(\widetilde{\lambda}\ot\widetilde{\mu})\circ\Delta\quad(\lambda,\mu\in A^\ast),\]
such that $A^\ast$ becomes a dual Banach algebra \emph{i.e.} $``\star"$ is separately weak*-continuous, and $A$ becomes a $A^\ast$-$A^\ast$ bi-module:  
\[
a\star \lambda:=(\widetilde{\lambda}\ot\id)(\Delta(a))\in A~;~\lambda\star a:=(\id\ot\widetilde{\lambda})(\Delta(a))\in A\quad(\forall~a\in A,\lambda\in A^\ast). 
\]
For $a\in A$ and $\lambda\in A^\ast$, we define the functionals $\lambda a := \lambda (a -)$ and $a\lambda := \lambda (-a)$. 

\bdfn\label{definition of Haar state}
A state $h\in A^\ast$ is an invariant mean if \[(\id\ot\widetilde{h})(\Delta(a))=h(a)1=(\widetilde{h}\ot\id)(\Delta(a)) \quad(\forall~a\in A).\]
\edfn

In this section, we will show that a compact semitopological quantum semigroup satisfying the following assumptions is a compact quantum group in the sense of \cite{Woronowicz}. Our proofs are the semitopological counterpart of the proof of Theorem 3.2 in \cite{weak cacellation}.

\begin{ass}\label{ass 1}
There exists an invariant mean $h\in A^\ast$ on $\mathbb{S}$.
\end{ass}

\begin{ass}\label{ass 2}
 The sets $\linear\{a\star hb:~a,b\in A\}$ and $\linear\{ha\star b:~a,b\in A\}$ are norm dense in $A$.
\end{ass}

In Sections \ref{QuEllis} and \ref{QuConv} we will study concrete, natural situations where these assumptions are satisfied.

\subsection{Preliminary results}

In this subsection, we consider a compact semitopological quantum semigroup $\mathbb{S}:=(A,\Delta)$ satisfying the assumptions \ref{ass 1} and \ref{ass 2}. 

The following result may be well-known, but we include a proof for the sake of completeness.

\begin{lmma}\label{Vaes-lemma}
Let $H,K$ be Hilbert spaces and $U,V\in B(H\ot K)$. Suppose $\{e_i\}_{i\in \cli}$ is an orthonormal basis for $H$ and $\xi\in H$. Let $p_i:=(\omega_{e_i,\xi}\ot\iota)(U)$ and $q_i:=(\omega_{\xi,e_i}\ot\iota)(V)$. Then for any $L\subset \cli$ we have 
\[
\|\sum_{i\in L}p_iq_i\|^2\leq\|\sum_{i\in L}p_ip_i^*\|\|\sum_{i\in L}q_i^*q_i\|. 
\]
\end{lmma}
\begin{proof}
For $i\in \cli$ let $P_i$ be the rank one projection of $H$ onto $\IC e_i$. Then $\sum_{i\in \cli}P_i=\id_{H}$ where the sum converges in the SOT* topology of $B(H)$. It follows that the sum $\sum_{i\in L}P_i$ converges in the SOT* topology as well. From this we can conclude that the series $\sum_{i\in L}(\omega_{\xi,\xi}\ot\id)(U(P_i\ot 1)V)=\sum_{i\in L}p_iq_i$ converges in the SOT* topology of $B(K)$. Similar arguments hold for the series $\sum_{i\in L}p_ip_i^*$ and $\sum_{i\in L}q_i^*q_i$, so that the RHS and LHS of the above inequality are finite.

Let us estimate the norm of the operator $\sum_{i\in L}p_iq_i$. Let $u,v\in K$ such that $\|u\|\leq1$ and $\|v\|\leq1$. From the above discussions it follows that the series $\sum_{i\in L}\innerl (p_iq_i)u,v\innerr$ consisting of scalars is convergent. Moreover, we have that 
\begin{equation*}
\begin{split}
\sum_{i\in L}|\innerl p_iq_iu,v\innerr|&\leq\sum_{i\in L}\|q_iu\|\|p_i^*v\|\\
&\leq(\sum_{i\in L}\|q_iu\|^2)^{\frac{1}{2}}(\sum_{i\in L}\|p_i^*v\|^2)^{\frac{1}{2}}\\
&\leq(\|\sum_{i\in L}p_ip_i^*\|)^{\frac{1}{2}}(\|\sum_{i\in L}q_i^*q_i\|)^{\frac{1}{2}}.
\end{split}
\end{equation*}
Thus we have 
\begin{equation*}
\begin{split}
|\innerl\sum_{i\in L}p_iq_i u, v\innerr|&=|\sum_{i\in L}\innerl p_iq_i u, v\innerr|\\
&\leq\sum_{i\in L}|\innerl p_iq_iu,v\innerr|\\
&\leq(\|\sum_{i\in L}p_ip_i^*\|)^{\frac{1}{2}}(\|\sum_{i\in L}q_i^*q_i\|)^{\frac{1}{2}}
\end{split}
\end{equation*}
from which the result follows.
\end{proof}
We will also be using the following:
\begin{ppsn}\label{Kus-Vaes-lemma}(Lemma A.3 in \cite{vaes})
Let C be a unital C*-algebra and $(x_\alpha)_\alpha\subset C$ be an increasing net of positive elements such that there exists a positive element $x\in C$ so that 
$\omega(x)=\text{sup}\{\omega(x_\alpha)\}$ for all states $\omega$. Then $x_\alpha\longrightarrow x$ in norm.
\end{ppsn}

\begin{lmma}\label{norm closure contains injective tensor product}
Let $\ot$ denote the injective tensor product of C*-algebras. Then we have 
\begin{enumerate}
\item
$A\ot A\subset \overline{(A\ot1)\Delta(A)}^{\|\cdot\|_{A}}$,
\item
$A\ot A\subset\overline{(1\ot A)\Delta(A)}^{\|\cdot\|_{A}}$;
\end{enumerate}
where $\|\cdot\|_A$ is the norm in $A$.
\end{lmma}

\begin{proof}
For $a,b\in A$  let $U:=\Delta(b)\ot 1$ and $V:=(\tilde{\Delta}\ot\id)(\Delta (a))$. Let $h$ denote the invariant state of $\mathbb{S}$. Considering $A\subset B(\clk)$ where $\clk$ is the universal Hilbert space of $A$, we see that there exists $\xi\in\clk$ such that $h:=\omega_{\xi,\xi}$. We have
\begin{equation*}
\begin{split}
(\widetilde{h}\ot\id\ot\id)(UV)&=(\widetilde{h}\ot\id\ot\id)\big((\Delta\ot\id)((b\ot 1)(\Delta a))\big)\\
&=1\ot \big((\widetilde{h}\ot\id)((b\ot 1)(\Delta a))\big)\\
&=1\ot a\star h b.
\end{split} 
\end{equation*}
Note that $U,V\in A^{\ast\ast}\otol A^{\ast\ast}\otol A^{\ast\ast}$. Consequently we have
\[
(\widetilde{h}\ot\id\ot\id)(UV)=(\omega_{\xi,\xi}\ot\id\ot\id)(UV)=\sum_{i\in\cli}(\omega_{e_i,\xi}\ot\id\ot\id)(U)(\omega_{\xi,e_i}\ot \id\ot\id)(V), 
\]
where $(e_i)_{i\in\cli}$ is an orthonormal basis for $\clk$, the sum being convergent in the ultraweak topology of $A^{\ast \ast}\otol A^{\ast \ast}$. 

At this point we may observe that 
\[
(\omega_{e_i,\xi}\ot\id\ot\id)(U)=(\omega_{e_i,\xi}\ot\id)(\Delta (b))\ot 1.	
\]
The series $\sum_{i\in\cli}(\omega_{e_i,\xi}\ot\id)(\Delta (b))(\omega_{\xi,e_i}\ot\id)(\Delta (b^*))$ converges in the ultraweak topology of $A^{\ast\ast}$. Since $\mathbb{S}$ is a compact semitopological quantum semigroup, it follows that 
\[
a_F:=\sum_{i\in F}(\omega_{e_i,\xi}\ot\id)(\Delta (b))(\omega_{\xi,e_i}\ot\id)(\Delta (b^*))\in A\quad(F\subset\cli,~|F|<\infty). 
\]
Moreover, we see that the net $(a_F)_{F\subset\cli}$ is increasing and converges in the ultraweak topology of $A^{\ast\ast}$ to $(\omega_{\xi,\xi}\ot\id)(\Delta(bb^*))=h(bb^*)1\in A$. Thus, by Proposition \ref{Kus-Vaes-lemma} $(a_F)_{F\subset\cli}$ converges in the norm topology to $h(bb^*)1$. Thus, the net
\[
b_F:=\sum_{i\in F}(\omega_{e_i,\xi}\ot\id\ot\id)(U)(\omega_{\xi,e_i}\ot\id\ot\id)(U^*)\quad(F\subset\cli,~|F|<\infty) 
\]
converges in the norm topology. This observation and Lemma \ref{Vaes-lemma} state that the sum 
\[ \sum_{i\in\cli}(\omega_{e_i,\xi}\ot\id\ot\id)(U)(\omega_{\xi,e_i}\ot \id\ot\id)(V) \]
converges in the norm topology of $A^{\ast\ast}\otol A^{\ast\ast}$. 

For each $i\in\cli$, $(\omega_{e_i,\xi}\ot\id\ot\id)(U)\in A\ot1$ and
\begin{equation*}
\begin{split}
(\omega_{\xi,e_i}\ot\id\ot\id)(V)&=\big(\omega_{\xi,e_i}\ot\id\ot\id\big)(\tilde{\Delta}\ot\id)\Delta (a)\\
&=\big(\omega_{\xi,e_i}\ot\id\ot\id\big)(\id\ot\tilde{\Delta})\Delta (a)\\
&=\Delta\big((\omega_{\xi,e_i}\ot\id)(\Delta (a))\big)\in\Delta(A). 
\end{split}
\end{equation*}
Thus, $(\omega_{e_i,\xi}\ot\id\ot\id)(U)(\omega_{\xi,e_i}\ot \id\ot\id)(V)\in (A\ot1)\Delta(A)$ which shows that 
\[ 1\ot a\star h b\in \overline{(A\ot1)\Delta(A)}^{\|\cdot\|_A}. \]
Since $\mathbb{S}$ satisfies Assumption \ref{ass 2}, we have
\[ A\ot A\subset  \overline{(A\ot1)\Delta(A)}^{\|\cdot\|_A}. \]
We may repeat the same argument with $\mathbb{S}^{\text{op}}:=(A,\tau \circ\Delta)$, where $\tau :A^{\ast\ast}\otol A^{\ast\ast}\longrightarrow A^{\ast\ast}\otol A^{\ast\ast}$ is the flip, to conclude that 
\[ A\ot A\subset \overline{(1\ot A)\Delta(A)}^{\|\cdot\|_A}. \]This proves our claim.
\end{proof}
Let $\clh$ denote the GNS Hilbert space of $A^{\ast\ast}$ associated with $\widetilde{h}$ and $\xi_0$ denote the cyclic vector. We will adopt the convention that whenever we write $a\xi$ for $a\in A$ and $\xi\in \clh$, we mean that $A$ acts on $\clh$ via the GNS representation of $h$. Note that the set $\{a\xi_0:~a\in A\}$ is norm dense in $\clh$. As before, we consider $A\subset B(\clk)$ where $\clk$ is the universal Hilbert space of $A$.

We omit the proof of the following result, which is similar to the proof of Proposition 5.2 in \cite{vanDaele}, using here Lemma \ref{norm closure contains injective tensor product}.

\begin{lmma}\label{the left regular representation}
There exists a unitary operator $u$ on $\clk\ot\clh$ given by 
\[
u(\eta\ot a\xi_0):=\Delta(a)(\eta\ot\xi_0)\quad(\eta\in\clk,a\in A). 
\]
\end{lmma}
For a C*-algebra $\clb$, $M_l(\clb)$ and $M_r(\clb)$ will denote the set of left and right multipliers of $\clb$.

\begin{lmma}\label{u is a right multiplier}
The operator $u$ (resp. $u^*$) belongs to $M_r(A\ot B_0(\clh))$ (resp. $M_l(A\ot B_0(\clh))$).
\end{lmma}
\begin{proof}
We will prove that $u^*\in M_l(A\ot B_0(\clh))$. For $a\in A$ and $\xi_1\in \clh$, let $\theta_{\xi_1,a\xi_0}\in B_0(\clh)$ denote the rank one operator given by $\theta_{\xi_1,a\xi_0}(\xi):=\innerl\xi,\xi_1\innerr a\xi_0$. Consider the operator $u^*(1\ot\theta_{\xi_1,a\xi_0})\in B(\clk\ot\clh)$. Let $\eta\in\clk$ and $\xi\in\clh$. We have
\begin{equation*}
\begin{split}
\big(u^*(1\ot\theta_{\xi_1,a\xi_0})\big)(\eta\ot\xi)&=\big(u^*(1\ot a)(1\ot\theta_{\xi_1,\xi_0})\big)(\eta\ot\xi).
\end{split}
\end{equation*}
Now by Lemma \ref{norm closure contains injective tensor product}, we see that $1\ot a$ can be approximated in norm by elements of the form $\sum_{i=1}^k\Delta(b_i)(c_i\ot1)$ where $b_i,c_i\in A$ for $i=1,2,...,k$. By a direct computation we can verify that 
\[
\sum_i\big(\Delta(b_i)(c_i\ot 1)(1\ot\theta_{\xi_1,\xi_0})\big)(\eta\ot\xi)=u\big(\sum_i c_i\ot\theta_{\xi_1,b_i \xi_0}\big)(\eta\ot\xi)\quad(\eta\in\clk,\xi\in\clh),
\]
so that we have
\begin{equation*}
\begin{split}
\|u^*(1\ot\theta_{\xi_1,a\xi_0})-\sum_{i=1}^k c_i\ot\theta_{\xi_1,b_i\xi_0}\|&=
\|u^*(1\ot a)(1\ot\theta_{\xi_1,\xi_0})-u^*u\sum_i c_i\ot\theta_{\xi_1,b_i\xi_0}\|\\
&\leq\|\xi_1\|\|\xi_0\|\|(1\ot a)-\sum_{i=1}^k\Delta(b_i)(c_i\ot 1)\|\longrightarrow0.
\end{split}
\end{equation*}
Thus we ended up proving that $u^*(1\ot\theta_{\xi_1,a\xi_0})\in A\ot B_0(\clh)$. Since $\xi_1$ is arbitrary and $A$ acts non-degenerately on $\clh$, it follows that $u^*(1\ot x)\in A\ot B_0(\clh)$ for all finite rank operators $x\in B_0(\clh)$ which in turn implies that $u^*(1\ot x)\in A\ot B_0(\clh)$ for all $x\in B_0(\clh)$. This proves that $u^*\in M_l(A\ot B_0(\clh))$.
\end{proof}

\begin{lmma}\label{slices of u generates A}
For $\omega\in B(\clh)_\ast$, the set $\{(\id\ot\omega)(u):~\omega\in B(\clh)_\ast\}$ is norm dense in $A$.
\end{lmma}
\begin{proof}
Let us first show that $\{(\id\ot\omega)(u):~\omega\in B(\clh)_\ast\}\subset A$. Let $a,b\in A$ and $\eta_1,\eta_2\in \clk$. We have
\begin{equation*}
\begin{split}
\innerl(\id\ot\omega_{a\xi_0,b\xi_0})(u)\eta_1,\eta_2\innerr&=\innerl u(\eta_1\ot a\xi_0),\eta_2\ot b\xi_0\innerr\\
&=\innerl\Delta(a)(\eta_1\ot\xi_0),\eta_2\ot b\xi_0\innerr\\
&=\innerl(\id\ot h)(1\ot b^*)\Delta(a)\eta_1,\eta_2\innerr.
\end{split}
\end{equation*}
This proves that $(\id\ot\omega_{a\xi_0,b\xi_0})(u)=(\id\ot h)(1\ot b^*)\Delta(a)$ and since $\mathbb{S}$ is a compact semitopological quantum semigroup, it follows that $(\id\ot h)(1\ot b^*)\Delta(a)\in A$.

By Assumption \ref{ass 2}, elements of the form $(\id\ot \widetilde{h})(1\ot b^*)\Delta(a)$ are total in $A$ in the norm topology. This observation, coupled with the fact that $A$ acts non-degenerately on $\clh$ imply that the required set is norm dense in $A$.
\end{proof}

%
Moreover, we have:

\begin{lmma}\label{u belongs to Astarstar otol B(H)}
We have that $u\in A^{\ast\ast}\otol B(\clh)$. 
\end{lmma}
\begin{proof}
Define the following CB map from $B(\clh)_\ast\longrightarrow A^{\ast\ast}$ given by 
\[ \omega\longrightarrow(\id\ot\omega)(u)\quad(\omega\in B(\clh)_\ast). \]
From Lemma \ref{slices of u generates A} it follows that $(\id\ot\omega)(u)\in A\subset A^{\ast\ast}$. It now follows that $u\in A^{\ast\ast}\otol B(\clh)$. 
\end{proof}

\begin{lmma}\label{u is a unitary representation}
We have $(\widetilde{\Delta}\ot\id)(u)=u_{23}u_{13}$.
\end{lmma}
\begin{proof}
Let $b,c\in A$ and $\eta_1,\eta_2\in\clk$. We have
\begin{equation*}
\begin{split}
u_{23}(b\ot1\ot c)(\eta_1\ot\eta_2\ot \xi_0)&=u_{23}(b\eta_1\ot\eta_2\ot c\xi_0)\\
&=(\id\ot\Delta)(b\ot c)(\eta_1\ot\eta_2\ot\xi_0).
\end{split}
\end{equation*}
Now we can approximate $\Delta(a)$ in the weak operator topology of $B(\clk\ot\clk)$ by elements of the form $\sum_{i=1}^k b_i\ot c_i$, such that $b_i,c_i\in A$ for all $i$. Using this in the above equation we get:
\begin{equation*}
\begin{split}
u_{23}u_{13}(\eta_1\ot\eta_2\ot a\xi_0)=(\id\ot\Delta)\Delta(a)(\eta_1\ot\eta_2\ot\xi_0). 
\end{split}
\end{equation*}
On the other hand we have
\[ (\Delta\ot\id)(b\ot y)\big(\eta_1\ot\eta_2\ot a\xi_0\big)=\Delta(b)(\eta_1\ot\eta_2)\ot ya\xi_0.  \]
Since by Lemma \ref{u belongs to Astarstar otol B(H)} we have $u\in A^{\ast\ast}\otol B(\clh)$, again approximating $u$ by $\sum_{i=1}^k b_i\ot y_i$ where $b_i's\in A$ and $y_i's\in B(\clh)$ in the weak operator topology of $B(\clk\ot\clh)$, we may replace the left side of the last equation by 
\[ (\Delta\ot\id)(u)\big(\eta_1\ot\eta_2\ot a\xi_0\big).  \]
Let us consider the right side of the equality. Let $\eta_1^\prime,\eta_2^\prime\in\clk$ and $\lambda\in \clh$. We see that 
\[ \innerl\Delta(b)(\eta_1\ot\eta_2)\ot ya\xi_0,\eta_1^\prime\ot\eta_2^\prime\ot\lambda\innerr=\innerl(f\ot\id)(b\ot y)(a\xi_0),\lambda\innerr, \]
where $f:=\omega_{\eta_1\ot\eta_2,\eta^\prime_1\ot\eta^\prime_2}\circ\widetilde{\Delta}$. Replacing $f$ by a normal functional of the form $\omega_{\eta_1,\eta_2}$ for $\eta_1,\eta_2\in\clk$ and approximating $u$ by linear combinations of elements $b\ot y$ as in the left side, we have
\[ \innerl(\omega_{\eta_1,\eta_2}\ot\id)(u)(a\xi_0),\lambda\innerr=\innerl(\omega_{\eta_1,\eta_2}\ot\id)(\Delta(a))(\xi_0),\lambda\innerr. \]
Thus, for any normal functional $f\in B(\clk)_\ast$ we have 
\[ \innerl(f\ot\id)(u)(a\xi_0),\lambda\innerr=\innerl(f\ot\id)(\Delta(a))(\xi_0),\lambda\innerr. \]
Thus, taking $f=\omega_{\eta_1\ot\eta_2,\eta^\prime_1\ot\eta^\prime_2}\circ\widetilde{\Delta}$ we arrive at the equation 
\[ (\Delta\ot\id)(\Delta(a))(\eta_1\ot\eta_2\ot \xi_0)=(\Delta\ot\id)(u)(\eta_1\ot\eta_2\ot a\xi_0). \]
By coassociativity of $\Delta$, it follows that 
\[ (\Delta\ot\id)(u)(\eta_1\ot\eta_2\ot a\xi_0)=(\id\ot\Delta)\Delta(a)(\eta_1\ot\eta_2\ot\xi_0)=u_{23}u_{13}(\eta_1\ot\eta_2\ot a\xi_0).   \]
This proves the result.
\end{proof}

\begin{rmrk}
$(A,\Delta,u^*,\clh)$ is a C*-Eberlein-algebra in the sense of \cite{das-daws}. 
\end{rmrk}
We closely follow the techniques given in Section 6 of \cite{vanDaele}. We will borrow some standard notations from the representation theory of topological groups. Let $\text{Mor}(u)$ denote the set of all operators $x\in B(\clh)$ such that $u(1\ot x)=(1\ot x)u$, \emph{i.e.} operators $x\in B(\clh)$ which intertwines $u$ with itself.

\begin{lmma}\label{compact intertwiners exist}
We have $\text{Mor}(u)\cap B_0(\clh)\neq \{0\}$. 
\end{lmma}
\begin{proof}
Let $x\in B_0(\clh)$ and consider the operator $u^*(1\ot x)u\in A^{\ast\ast}\otol B(\clh)$. By Lemma \ref{u is a right multiplier}, it follows that $u^*(1\ot x)\in A\ot B_0(\clh)$. Again, applying Lemma \ref{u is a right multiplier} to $u^*(1\ot x)$ and considering $u$ as a right multiplier, we get $u^*(1\ot x)u\in A\ot B_0(\clh)$. Let $y:=(h\ot\id)(u^*(1\ot x)u)\in B_0(\clh)$. We have
\begin{equation*}
\begin{split}
(\Delta\ot\id)(u^*(1\ot x)u)&=u^*_{13}u^*_{23}(1\ot1\ot x)u_{23}u_{13}.
\end{split}
\end{equation*}
Applying $(\id\ot h\ot\id)$ to both sides of this equation and using the translation invariance of $h$ we arrive at 
\[ 1\ot y=u^*(1\ot y)u. \]
Since $u$ is unitary, we get $u(1\ot y)=(1\ot y)u$, and $y\in B_0(\clh)$. 

Suppose $y=0$ for all $x\in B_0(\clh)$. Let $\{x_\alpha\}_{\alpha\in \Lambda}\in B_0(\clh)$ be an approximate identity, converging to $1$ in the ultraweak topology of $B(\clh)$. It follows that \[0=y_\alpha:=(\widetilde{h}\ot\id)(u^*(1\ot x_\alpha)u)\longrightarrow (\widetilde{h}\ot\id)(1\ot 1)=1,\] a contradiction. Hence $y\neq0$ for some $x\in B_0(\clh)$, which proves the result.
\end{proof}

\begin{lmma}\label{Peter-Weyl decomposition}
There exists a set $\{e_\alpha:~\alpha\in \cli\}$ of mutually orthogonal finite-dimensional projections on $\clh$ with sum $1$ and satisfying
\[ u(1\ot e_\alpha)=(1\ot e_\alpha)u\quad(\forall~\alpha\in \cli). \]
\end{lmma}
\begin{proof}
Let $\clb:=\{y\in B_0(\clh):~u(1\ot y)=(1\ot y)u\}$. By Lemma \ref{compact intertwiners exist}, we have $\clb\neq\emptyset$. Moreover, $\clb$ is a norm closed subalgebra of $B(\clh)$. By Lemma \ref{the left regular representation}, $u$ is unitary which implies that $\clb$ is self-adjoint.  Thus $\clb$ is a C*-subalgebra of $B_0(\clh)$. 

By Lemma \ref{compact intertwiners exist}, $(h\ot\id)(u^*(1\ot x)u)\in \clb$ for all $x\in B_0(\clh)$. Let $(x_\lambda)_\lambda \subset B_0(\clh)$ be an increasing net of positive elements such that $x_\lambda\longrightarrow1$ in the ultraweak topology of $B(\clh)$. Then $y_\lambda:=(\widetilde{h}\ot\id)(u^*(1\ot x_\lambda)u)\longrightarrow (\widetilde{h}\ot\id)(u^*u)=1$. This implies that $\clb$ acts non-degenerately on $\clh$. Thus, we may select a maximal family of mutually orthogonal, minimal projections in $\clb$ say $\{e_\alpha:~\alpha\in\cli\}$ which are finite-dimensional as $\clb\subset B_0(\clh)$. Non-degeneracy of $\clb$ implies that $\oplus_{_{\alpha\in\cli}}e_\alpha=1$. This proves the assertion.
\end{proof}

\subsection{Main result}
We recall the definition of a compact quantum group from \cite{Woronowicz, vanDaele}.
\bdfn\label{definition of a compact quantum group}
A compact quantum group $\mathbb{G}:=(A,\Delta)$ consists of a unital C*-algebra $A$ and a unital *-homomorphism $\Delta: A\longrightarrow A\ot A$ such that
\begin{itemize}
 \item $\Delta$ is coassociative: $(\id\ot\Delta)\circ\Delta=(\Delta\ot\id)\circ\Delta$,
 \item $(1\ot A)\Delta(A)$ and $(A\ot 1)\Delta(A)$ are norm dense in $A\ot A$ (Woronowicz cancellation laws).
\end{itemize}
\edfn
The following result from \cite{Woronowicz,vanDaele} justifies the word ``quantum''.
\begin{ppsn}
Let $\mathbb{G}:=(A,\Delta)$ be a compact quantum group with $A$ commutative. Then $A=C(G)$ for a compact group $G$. 
\end{ppsn}

\bthm\label{S is a compact quantum group}
A compact semitopological quantum semigroup $\mathbb{S}:=(A,\Delta)$ satisfying the assumptions \ref{ass 1} and \ref{ass 2}, that is:
\begin{enumerate}
 \item there exists an invariant mean $h\in A^\ast$ on $\mathbb{S}$,
 \item the sets $\linear\{a\star hb:~a,b\in A\}$ and $\linear\{ha\star b:~a,b\in A\}$ are norm dense in $A$.
\end{enumerate} 
 is a compact quantum group.  
\ethm
\begin{proof}
From Theorem \ref{Peter-Weyl decomposition}, it follows that $u=\oplus_{\alpha\in\cli}u_\alpha$ where $u_\alpha:=u(1\ot e_\alpha)$. So each $u_\alpha\in A^{\ast\ast}\otol B(H_\alpha)$ where $H_\alpha:=e_\alpha(\clh)$ and also $\clh=\oplus_{\alpha\in\cli}H_\alpha$. Thus, from Theorem \ref{slices of u generates A} it follows that the linear span of the set 
$\{(\id\ot\omega)(u):~\omega\in B(H_\alpha)_\ast\,\alpha\in\cli\}$ is norm dense in $A$. Let $\text{dim}H_\alpha=m$ and $\{f_i\}_{i=1}^m$ be an orthonormal basis for $H_\alpha$. Note that we have $(\Delta\ot\id)(u_\alpha)=u_{\alpha 23}u_{\alpha 13}$. Taking $\xi,\eta\in H_\alpha$ we have
\[ \Delta((\id\ot\omega_{\xi,\eta})(u_\alpha))=\sum_{i=1}^m(\id\ot\omega_{f_i,\eta})(u_\alpha)\ot (\id\ot\omega_{\xi,f_i})(u_\alpha)\in A\ot A. 
\]
This implies that $\Delta(A)\subset A\ot A$. This observation and Lemma \ref{norm closure contains injective tensor product} now implies that the sets 
$\Delta(A)(A\ot1)$ and $\Delta(A)(1\ot A)$ are norm dense in $A\ot A$. Thus $\mathbb{S}:=(A,\Delta)$ is a compact quantum group.\end{proof}
\begin{rmrk}\label{quantum semigroups are semitopological}
Note that in Definition \ref{definition of compact semitopological quantum semigroup} if we take $\Delta:A\longrightarrow A\ot A$, then the resulting object $\mathbb{S}:=(A,\Delta)$ is a compact quantum (topological) semigroup as defined in \cite{weak cacellation}. Thus, Theorem \ref{S is a compact quantum group} generalizes Theorem 3.2 in \cite{weak cacellation}. 
\end{rmrk}
The next corollary gives a situation when a compact semitopological quantum semigroup satisfies Assumption \ref{ass 1} but need not satisfy Assumption \ref{ass 2}. This shows that in general there is no relation between the two assumptions. 

Let us recall a few facts on quantum Bohr compactification and quantum Eberlein compactification of quantum groups (see \cite{soltan,das-daws}). Let $\qg$ be a locally compact quantum group and $(C_0^u(\qg),\Delta_u)$ be its universal version. Let $(\ap(C_0^u(\qg)),\Delta)$ denote its quantum Bohr compactifiction (see \cite[Proposition 2.13]{soltan}) and $E(\qg):=(\cle(\qg),\Delta_\qg,U_\qg,H_\qg)$ be its quantum Eberlein compactification (see \cite[Theorem 5.2]{das-daws}). By Theorem \ref{answer to Matt's question}, we see that $(\cle(\qg),\Delta_\qg)$ is a compact semitopological quantum semigroup. By construction $\ap(C_0^u(\qg))\subset\cle(\qg)\subset M(C_0^u(\qg))$, $\Delta=\overline{\Delta}_u|_{_{\ap(C_0^u(\qg))}}$ and $\Delta_\qg=\overline{\Delta}_u|_{_{\cle(\qg)}}$, where $\overline{\Delta}_u$ denotes the strict extension of $\Delta_u$ to $M(C_0^u(\qg))$.

\begin{crlre}\label{examples of CSQS revisited}
The following facts are true:
\begin{itemize}
\item[(a)]
The compact semitopological quantum semigroup $(\cle(\qg),\Delta_\qg)$ satisfies Assumption \ref{ass 1}.
\item[(b)]
Let $(\ap(C_0^u(\qg)),\Delta)$ be reduced. Then $(\cle(\qg),\Delta_\qg)$ satisfies Assumption \ref{ass 2} if and only if $\qg$ is a coamenable compact quantum group.
\end{itemize}
\end{crlre}
\begin{proof}
From \cite[Theorem 4.6]{das-daws} it follows that $(\cle(\qg),\Delta_\qg)$ satisfies Assumption \ref{ass 1}, which proves (a). 

Let $h$ denote the invariant mean on $(\cle(\qg),\Delta_\qg)$, and assume it satisfies Assumption \ref{ass 2}. Suppose $\qg$ is non-compact. Then by Theorem \ref{S is a compact quantum group} it follows that $(\cle(\qg),\Delta_\qg)$ is a compact quantum group. From \cite[Proposition 2.13]{soltan} and the discussions in \cite[Section 3]{soltan} it follows that $\cle(\qg)=\ap(C_0^u(\qg))$. From the remarks after \cite[Theorem 5.3]{das-daws} it follows that $C_0^u(\qg)\subset \cle(\qg)=\ap(C_0^u(\qg))$. From \cite[Proposition 5.4]{das-daws} it follows that $h$ annihilates $C_0^u(\qg)$. Since $(\ap(C_0^u(\qg)),\Delta)$ is reduced, this implies that $C_0^u(\qg)=\{0\}$. This leads to a contradiction. Hence $\qg$ must be compact. Coamenability of $\qg$ now follows from the fact that $\ap(C^u(\qg))=C^u(\qg)$ and $(\ap(C^u(\qg)),\Delta)$ is reduced.

Suppose $\qg$ is coamenable and compact. Compactness of $\qg$ implies (by \cite[Theorem 5.1]{das-daws}) that $\cle(\qg)=\ap(C^u(\qg))=C^u(\qg)$. Thus $(\cle(\qg),\Delta_\qg)$ becomes a compact quantum group and it satisfies Assumption \ref{ass 2}. This proves (b). 
\end{proof}
\begin{rmrk}\label{remarks about when as1 and or as2}
Corollary \ref{examples of CSQS revisited} implies the following: 
\begin{itemize}
 \item[(a)]
Let $G$ be a locally compact, non-compact group. Suppose $G^\cle$ is the Eberlein compactification of $G$. Then the commutative compact semitopological quantum semigroup $(C(G^\cle),\Delta)$ (see Section \ref{subsection1} for the definition of $\Delta$) satisfies Assumption \ref{ass 1} but does not satisfy Assumption \ref{ass 2}.
\item[(b)]
Suppose $G$ is a locally compact, non-compact, non-discrete and non-abelian group, such that $G_d$ ($G$ with discrete topology) is amenable. Let $\widehat{G}:=(C^*(G),\Delta_{\widehat{G}})$ denote the universal dual quantum group of $G$. From the discussions in \cite[Section 8]{das-daws} it follows that the Eberlein compactification of $\widehat{G}$ namely $\cle(\widehat{G})$ is the closure of the image of the measure algebra $M(G)=C_0(G)^*$ inside $M(C^*(G))$. Moreover, letting $\Delta^\prime:=\overline{\Delta}_{\widehat{G}}|_{_{\cle(\widehat{G})}}$ it follows that $(\cle(\widehat{G}),\Delta^\prime)$ is a noncommutative compact semitopological quantum semigroup. The fact that $G_d$ is amenable implies that $(\ap(\widehat{G}),\Delta)$ (the quantum Bohr compactification of the universal dual quantum group $\widehat{G}$) is reduced (see \cite[Proposition 4.3]{soltan}). Thus by Corollary \ref{examples of CSQS revisited} it follows that the noncommmutative compact semitopological quantum semigroup $(\cle(\widehat{G}),\Delta^\prime)$ satisfies Assumption \ref{ass 1} but does not satisfy Assumption \ref{ass 2}.
\end{itemize}
\end{rmrk}
\begin{rmrk}
The compact semitopological quantum semigroup $(\cle(\widehat{G}),\Delta^\prime)$ appearing in Remark \ref{remarks about when as1 and or as2} (b) is cocommutative i.e. $\Delta^\prime=\tau\circ\Delta^\prime$, where $\tau$ is the tensor flip. It is a well-known fact (see \cite[Theorem 4.2.4]{enock}) that all cocommutative locally compact quantum groups are precisely group C*-algebras of locally compact groups. It is tempting to conjecture that all cocommutative compact semitopological quantum semigroups also arise from group C*-algebras of locally compact groups, somewhat as in Remark \ref{remarks about when as1 and or as2} (b). However at this stage, it is not clear to us whether such characterizations are possible in general.
\end{rmrk}

In the following sections, we consider situations where we can apply Theorem \ref{S is a compact quantum group}. Our motivations are results by Ellis in \cite{Ellis} and by Mukherjea and Tserpes in \cite{mukherjea}. In particular, we provide new proofs of these results in the compact case.

%
%

\section{Quantum Ellis joint continuity theorem}\label{QuEllis}

In \cite{Ellis}, Ellis showed that a compact semitopological semigroup which is algebraically a group is a compact group. A simplified proof of this was given in \cite{glicksberg}. The key point in this proof was to show that a group which is also a compact semitopological semigroup, always admits a faithful invariant mean. In fact the compact case plays an important role in the theory of weakly almost periodic compactification, in particular, in the structure theory of the kernel of a semigroup.

We will prove an analogous result for compact semitopological quantum semigroup, which in particular will give a new (C*-algebraic) proof of Ellis Theorem (Corollary \ref{Ellis thm}). However, for that we first need a noncommutative analogue of the condition ``algebraically a group". This is discussed in the following paragraph.

In what follows, $S$ will denote a compact semitopological semigroup and $\mathbb{S}$ will denote a compact semitopological quantum semigroup. 

\subsection{A necessary and sufficient condition on $C(S)$ for $S$ to be algebraically a group}
Our aim here is to identify a necessary and sufficient condition on $C(S)$ which implies that $S$ is algebraically a group. 
\bdfn\label{weak cancellation}
A compact semitopological quantum semigroup $\mathbb{S}:=(A,\Delta)$ is said to have 
\emph{weak cancellation laws} if it satisfies:
\[
\overline{\linear\{a\star\omega b:~a,b\in A\}}^{\|\cdot\|_{A}}=A=
\overline{\linear\{\omega b\star a:~a,b\in A\}}^{\|\cdot\|_{A}}
\]
for every state $\omega\in A^\ast$.
\edfn

\begin{rmrk}
 The weak cancellation laws in Definition \ref{weak cancellation} are inspired by \cite{weak cacellation}.
\end{rmrk}

In the classical case, we have the following.

\begin{thm}\label{A necessary condition for being algebraically a group}
 Let $S$ be a compact semitopological semigroup. The following are equivalent:
\begin{enumerate}
 \item\label{alg} $S$ is algebraically a group,
 \item\label{we} $(C(S), \Delta)$ has weak cancellation laws.
\end{enumerate}
\end{thm}

\begin{proof}
Assume that $S$ is algebraically a group. Note that the pure states on $C(S)$ are the evaluation maps $ev_y$, $y \in S$. Fix an element $y\in S$. It follows that the map $R_y:S\longrightarrow S$ given by $x\mapsto xy$ is a homeomorphism. Thus $R_y^\ast:C(S)\longrightarrow C(S)$ given by $R_y^\ast(f)(s):=f(R_y(s))$ is a C*-algebra isomorphism. It is easy to check that $R_y^\ast(f)=(\id \ot \widetilde{ev_y})(\Delta(f))$. Thus $C(S)=\text{Ran}~R_y^\ast$ = $\overline{\linear\{ev_y\star f:~f\in C(S)\}}^{\|\cdot\|_\infty}$. Similarly considering the map $L_y:S\longrightarrow S$ given by $x\mapsto yx$, we have $C(S)=\text{Ran}~L_y^\ast$ = $\overline{\linear\{f\star ev_y:~f\in C(S)\}}^{\|\cdot\|_\infty}$.

Let $\omega \in C(S)^*$ be a state. Since $C(S)$ is commutative, it follow that the set $\cli:=\{f\in C(S):~\omega(f^*f)=0\}$ is a 2-sided ideal in $C(S)$. Let $(\pi,\xi, H)$ denote the GNS triple associated with the state $\omega$ and let $A_r:=\pi(C(S))$. It is easy to see that $\cli$ being a 2-sided ideal, ker$\pi$ $=$ $\cli$. So the functional $\omega_r\in A_r^*$ defined by $\omega_r(\pi(f)):=\omega(f)$ for $f\in C(S)$ is a well-defined faithful state on $A_r$. Put $L:=\overline{\text{Lin}\{\omega f\star g:~f,g\in C(S)\}}^{\|\cdot\|_\infty}$, and let $\mu\in C(S)^*$ be such that $\mu(L)=0$. In particular, we have, for all  $f,g\in C(S)$, \[\mu(\omega f\star g)=(\widetilde{\mu}\ot \widetilde{\omega f})(\Delta(g))=0.\]	
Rewriting the last equation in terms of $\omega_r$ we have, for all $x,y\in C(S)$,
\begin{equation*}
\begin{split}
0&=(\widetilde{\mu}\ot \widetilde{\omega f})(\Delta(g))\\
&=(\widetilde{\mu}\ot \widetilde{\omega_r\pi(f)})\big((\id\ot\tilde{\pi})(\Delta(g))\big)\\
&=(\widetilde{\mu}\ot \widetilde{\omega_r}\pi(f))\big((\id\ot\widetilde{\pi})(\Delta(g))\big)\\
&=\omega_r\big(\pi(f)((\widetilde{\mu}\ot\widetilde{\pi})(\Delta(g))\big).
\end{split}
\end{equation*}
The fact that $\omega_r$ is faithful on $A_r$ and Cauchy-Schwarz inequality imply that, for all $g\in C(S)$,
\begin{equation}\label{the final step in proving weak cancellation}
(\widetilde{\mu}\ot\widetilde{\pi})(\Delta(g))=0). 
\end{equation}
Since $A_r$ is commutative, there exists a nonzero multiplicative functional $\Lambda\in A_r^*$. Thus $\Lambda\circ\pi:C(S)\longrightarrow\IC$ is a non-zero multiplicative, bounded linear functional. Thus there exists $s\in S$ such that $\Lambda\circ\pi=ev_{s}$. Applying $(\id\ot\Lambda)$ to equation \eqref{the final step in proving weak cancellation} we have, for all $g\in C(S)$,
\[
\mu\big((\id\ot\widetilde{\Lambda\circ\pi})(\Delta(g))\big)=\mu\big((\id\ot \widetilde{ev_s})(\Delta(g))\big)=\mu(ev_s\star f)=0.
\]
The fact that $\text{Ran}~R_s^*=C(S)$ implies that $\mu(f)=0$ for all $f\in C(S)$, and we must have $C(S)=\overline{\text{Lin}\{\omega f\star g:~f,g\in C(S)\}}^{\|\cdot\|_\infty}$. Similarly we can show that $\overline{\text{Lin}\{g\star \omega f:~f,g\in C(S)\}}^{\|\cdot\|_\infty}=C(S)$. Since $\omega\in C(S)^*$ was arbitrary it follows that $(C(S),\Delta)$ has weak cancellation laws. 

Let us introduce the kernel of $S$, denoted by $K(S)$, defined as the intersection of all two sided ideals of $S$. By Theorem 2.1 in \cite{burckel}, $S$ has minimal left and right ideals. Moreover, each minimal left or right ideal is closed. This fact coupled with Theorem 2.2 in \cite{burckel} imply that $K(S)\neq\emptyset$.

We first show that $S$ has right and left cancellation laws. For $p,q,r\in S$ let $pq=pr$. By the hypothesis $C(S)$ satisfies
\[
\overline{\{f\star ev_pg:~f,g\in C(S)\}}^{\|\cdot\|_\infty}=C(S).
\]
We have
\[
(f\star ev_pg)(q)=g(p)f(pq) 
\]
and 
\[
(f\star ev_p)(r)=g(p)f(pr).  
\]
The equality $pq=pr$ implies that $(f\star ev_pg)(q)=(f\star ev_pg)(r) $ for all $f,g\in C(S)$. The hypothesis that $\linear\{f\star ev_pg:~f,g\in C(S)\}$ is a norm dense subset of $C(S)$ yields $f(q)=f(r)$ for all $f\in C(S)$ which proves that $q=r$. Thus $S$ has left cancellation. 

Similarly using the other density condition in the hypothesis we can prove that $S$ has right cancellation.

So $S$ is a compact semigroup with right and left cancellations. We complete the proof by showing that $S$ has an identity and every element in $S$ has an inverse. 

For any $x\in K(S)$, since $xK(S)\subset K(S)$ and $xK(S)$ is a closed ideal in $S$, we have $xK(S)=K(S)$. So there exists $e\in K(S)$ such that $xe=x$. Multiplying to the right by $y\in S$ and using the fact that $S$ has left cancellation, we get that $ey=y$ for all $y\in S$. Then multiplying the last equation by any element $a\in S$ from the left and using the fact that $S$ has right cancellation, we get $ae=a$ for all $a\in S$. Thus $e$ is the identity of $S$ and in particular $K(S)=S$.

Let $s\in S$. As before we have $sK(S)=sS=K(S)=S$, so that there exists $p\in S$ such that $sp=e$. Thus $s$ has a left inverse. Similarly we can argue that $s$ has a right inverse. Since $s\in S$ was arbitrary, it follows that every element of $S$ has an inverse. Thus $S$ is algebraically a group.\end{proof}

\subsection{Compact semitopological quantum semigroup with weak cancellation laws}\label{here we prove quantum Ellis theorem}

As a consequence, we obtain a generalization of Ellis joint continuity theorem. But before that let us make one crucial observation concerning the existence of the Haar state.

\bthm\label{Ryll-Nardzewski in the compact case}
A compact semitopological quantum semigroup $\mathbb{S}$ with weak cancellation laws admits a unique invariant mean.
\ethm

\begin{proof}
We briefly remark on the arguments needed for this proof, which are similar to the arguments in \cite{weak cacellation}. Suppose $\mathbb{S}$ satisfies weak cancellation laws. We may repeat the exact arguments in the proofs of \cite[Lemma 2.1, Lemma 2.3, Theorem 2.4 and Theorem 2.5]{weak cacellation} to obtain the existence of a unique bi-invariant state $h$ on $\mathbb{S}$.\end{proof}


\bthm\label{Quantum Ellis theorem}({\bf Quantum Ellis joint continuity theorem}):
A compact semitopological quantum semigroup $\mathbb{S}:=(A,\Delta)$ is a compact quantum group if and only if it has the weak cancellation laws. 
\ethm
 
\begin{proof}
By Theorem \ref{Ryll-Nardzewski in the compact case}, $\mathbb{S}$ has a unique invariant mean, say $h$. By the hypothesis the two sets $\{ha\star b:~a,b\in A\}$ and $\{a\star hb:~a,b\in A\}$ are total in $A$ in norm. Thus by Theorem \ref{S is a compact quantum group}, $\mathbb{S}$ is a compact quantum group.  The converse easily follows by observing that Woronowicz cancellation laws imply weak cancellation laws.\end{proof}

Specializing to the commutative case, we have a new proof of Ellis joint continuity theorem.
\begin{crlre}\label{Ellis thm}({\bf Ellis joint continuity theorem}):
Let $S$ be a compact semitopological semigroup which is algebraically a group. Then $S$ is a compact group. 
\end{crlre}

\begin{proof}
On $C(S)$ define the map $\Delta:C(S)\longrightarrow C(S)^{\ast\ast}\otol C(S)^{\ast\ast}$ by $\Delta(f)(s,t):=f(st)$. From the discussions in Subsection \ref{subsection1} it follows that $(C(S),\Delta)$ is a compact semitopological quantum semigroup. It follows from Theorem \ref{A necessary condition for being algebraically a group} that $(C(S),\Delta)$ has weak cancellation laws. Hence by Theorem \ref{Quantum Ellis theorem} $(C(S),\Delta)$ is a compact quantum group which implies that $S$ is a compact group.\end{proof}

\subsection*{A comparison with the classical proof of Ellis theorem in compact case}
In \cite{glicksberg} the authors gave a proof of Ellis theorem for compact semitopological semigroups by using tools from the theory of weakly almost periodic functions on topological groups. We give a brief account of their proof, since it is the closest in spirit to our techniques.

Let $S$ be a locally compact semitopological semigroup, which is algebraically a group. Suppose $WAP(S)$ and $AP(S)$ denote respectively the algebra of weakly almost periodic and almost periodic functions on $S$. It is a highly non-trivial fact in the theory of weakly almost periodic functions that $WAP(S)$ admits a unique invariant mean $m$, such that for $f\in AP(S)$ with $f>0$, $m(f)>0$. The existence of such a mean can be proven using Ryll-Nardzewski fixed point theorem \cite[Corollary 1.26]{burckel}. If in addition $S$ is compact, it follows from the theory of weakly almost periodic functions that $WAP(S)=AP(S)=C(S)$, which implies that the invariant mean on $C(S)$ is faithful. 

Let $e$ be the identity of $S$ and assume that the map $S\times S\ni(\sigma,\tau)\mapsto \sigma^{-1}\tau\in S$ is not continuous. This means that there exists a neighbourhood $W\ni e$ such that for all neighbourhoods $U\ni e$, there exists $(\sigma_U,\tau_U)\in U\times U$ such that $\sigma_U^{-1}\tau_U\notin W$. Let $N$ be the directed set of neighbourhoods of $e$ with ordering by inclusion. Compactness of $S$ allows us to get a net $\{\eta_U:=\sigma_U^{-1}\tau_U\}_U$ such that $\lim_N\eta_U =\eta\neq e$. Let $f\in C(S)$ be such that $f(e)\neq f(\eta)$. For $s\in S$, denoting the right translation of $f$ by $s$ as $R_sf$, it follows that $f\neq R_\eta f$. Faithfulness of the mean $m$ implies that $\innerl|f-R_\eta f|,m\innerr>0$. Using this and some more facts about weakly almost periodic functions, one arrives at $\lim_N\innerl |R_{\eta_U}f-R_{\eta}f|,m\innerr>0$. The proof is now completed by showing that the net of real numbers $\{\innerl |R_{\eta_U}f-R_{\eta}f|,m\innerr\}_N$ has $0$ as a limit point which leads to a contradiction. 

Theorem \ref{Ryll-Nardzewski in the compact case} is a quantum generalization of the result about the existence of invariant mean on $C(S)$, where $S$ is a compact semitopological group. However, unlike the classical case, this mean can have non-trivial kernel. For example consider a non-coamenable universal compact quantum group, where the Haar state will always have non-trivial kernel. So in particular, Corollary \ref{Ellis thm} gives a proof of Ellis theorem, where one does not require the faithfulness of the invariant mean, and hence is fundamentally different from the above proof.

It is worthwhile to mention that there are proofs of Ellis theorem in the compact case, using other sophisticated tools. For example in \cite[Theorem 2.2]{megrelishvili}, the authors gave a proof of Ellis theorem in the compact case using the theory of enveloping semigroups. However, such methods are far from our considerations. 

\section{Converse Haar's theorem for compact semitopological quantum semigroups}\label{QuConv}

The converse of Haar's theorem states that a complete separable metric group which admits a locally finite non-zero right (left) - invariant postive measure is a locally compact group, with the invariant measure being the right (left) Haar measure of the group. In \cite{mukherjea} it was shown (see Theorem 1. (b)) that a locally compact semitopological semigroup admitting an invariant mean with full support is a compact group. We will prove a similar result for a compact semitopological quantum semigroup. 

\bdfn
Let $\mathbb{S}:=(A,\Delta)$ be a compact semitopological quantum semigroup. A bounded counit for $\mathbb{S}$ is a unital $*$-homomorphism $\varepsilon : A \to \IC$ such that $(\widetilde{\varepsilon} \otimes \id)\Delta= \id= (\id \otimes \widetilde{\varepsilon})\Delta$.
\edfn

\bthm\label{faithful invariant mean gives CQG}
Let $\mathbb{S}:=(A,\Delta)$ be a compact semitopological quantum semigroup, admitting a faithful invariant state $h$ and a bounded counit $\varepsilon$. Then $\mathbb{S}$ is a coamenable compact quantum group.
\ethm

\begin{proof}
%
If we show that $\mathbb{S}$ is a compact quantum group, its coamenability will follow from Theorem 2.2 in \cite{coamenability}. By hypothesis, $\mathbb S$ satisfies Assumption \ref{ass 1} of Section \ref{main results}. We will show that it satisfies Assumption \ref{ass 2}. Let $L$ be the closed linear span of elements of the form $a\star hb$, with $a,b \in A$, and $R$ be the closed linear span of elements $ha \star b$. Let us show that $L=R=A$. We will only show that $L=A$, the proof for $R=A$ being identical. Let $\omega \in A^*$ be a non-zero linear functional on $A$ that vanishes on $L$. In particular, we have for all $a,b \in A$ $$0=\omega (a\star h b)=h\big(b(\id \otimes \widetilde\omega) \Delta(a)\big),$$ and since $h$ is faithful, this implies that $(\id \otimes \widetilde\omega) \Delta(a)=0$. Applying the counit we get $$0=\varepsilon \big((\id \otimes \widetilde\omega) \Delta(a)\big)=\omega \big((\widetilde\varepsilon \otimes \id)\Delta(a)\big)=\omega(a).$$ This implies that $\omega=0$, which proves that $L=A$ and consequently $\mathbb{S}$ satisfies Assumption \ref{ass 2}. Therefore it follows from Theorem \ref{S is a compact quantum group} that $\mathbb{S}$ is a compact quantum group.
\end{proof}
Restricting to the commutative case, we have
\begin{crlre}[Theorem 1.(b) in \cite{mukherjea}]
A compact semitopological semigroup with identity admitting a non-zero invariant mean with full support is a compact group. 
\end{crlre}
By virtue of Remark \ref{quantum semigroups are semitopological}, we also have the following as a special case (Theorem 4.2 in \cite{coamenability}):
\begin{crlre}
Let $\mathbb{S}:=(A,\Delta)$ be a compact (topological) quantum semigroup with a bounded counit, admitting a faithful invariant mean. Then $\mathbb{S}$ is a coamenable compact quantum group.  
\end{crlre}

\end{document}